\documentclass{article}

\usepackage{Preamble}
\usepackage[style=ieee,citestyle=numeric-comp]{biblatex}

\hypersetup{
  pdftitle={Stability of skill-based queues with deterministic waiting time thresholds},
  pdfauthor={Sanne van Kempen, Elene Anton, Fiona Sloothaak},
  pdfkeywords={Skill-based queues, Stability, Threshold policy, Piecewise deterministic Markov process}
}

\begin{document}
\title{Stability of skill-based queues with deterministic waiting time thresholds}
\author{Sanne van Kempen, Elene Anton, Fiona Sloothaak}
\date{}
\maketitle

\begin{abstract}
    We consider skill-based queues where service is First--Come--First--Served, but some compatibility lines are only available once a customer's waiting time exceeds a threshold.
    For this model, we establish necessary and sufficient stability conditions.
    We find that the thresholds do not affect stability; an intuitive result that has not been proved in the literature up to this point. 
    While necessity follows from a coupling argument, the sufficiency proof is more challenging.
    It follows by formulating the First--in--Line waiting time process as a Piecewise Deterministic Markov Process (PDMP)  and constructing a Lyapunov function that accounts for the threshold structure. 
    We show that this Lyapunov function satisfies the boundary condition of the PDMP framework, and that its generator has negative drift when waiting times are large. 
    Finally, we show that bounded waiting-time sets are petite, so that the Foster--Lyapunov criterion applies. 
\end{abstract}

\section{Introduction}
In skill-based service systems, there is an inherent trade-off between matching quality and delay.
Routing customers to their preferred server can improve service quality, but may also increase waiting times when the preferred servers are busy.
Threshold policies manage this trade-off by allowing additional service capacity once queues become congested.
We consider a service policy with \gls{FCFS} service and waiting time thresholds: a customer may  be served by certain compatible servers only if its waiting time exceeds a deterministic threshold. 
A schematic example of such a skill-based queue is shown in \Cref{fig:intro_qsystem_example}.

We prove the necessary and sufficient stability conditions for this threshold policy.
The most challenging aspect regards proving sufficiency, which we do by modeling the \gls{FIL} waiting time process as a \gls{PDMP}~\cite{Davis1984} and constructing a suitable Lyapunov function.
This framework fits our model naturally: the \gls{FIL} waiting time increases linearly over time and jumps down when a customer enters service, which is caused by either a departure, or when the customer's waiting time reaches the threshold. 
We find that the thresholds do not affect stability: the stability region agrees with the corresponding region under \gls{FCFS} service without thresholds~\cite{Adan2014}.
This result is not surprising:  when the system operates close to its capacity, queues build up and waiting times become large, so that thresholds no longer restrict customer routing.
However, the formal proof is non-trivial and the \gls{PDMP} framework we develop may be of independent interest for related threshold models.
Our result relies on an assumption that each server has at most one threshold connection. 
We discuss in \Cref{sec:assumption} what is needed to relax this assumption.

In the literature, there are several methods for analyzing stability in queues. 
A first method is to characterize the stationary distribution of the number of customers in the system and verify that its normalization constant is finite.
For skill-based queues under \gls{FCFS} service without waiting time thresholds, product form stationary distributions of the queue length process are known \cite{Adan2009,Visschers2012,Adan2012}. 
However, in our model the number of customers is not itself a Markov chain, since the service policy  is waiting time dependent. 
Hence, analyzing the stationary distribution of the number of customers is not straightforward. 

A second method is a coupling technique, where an upper bound system with known stability region is shown to dominate the original system on every sample path. 
For our model, a natural upper bound system is the following: each customer waits its largest threshold deterministically, and is then eligible at all its compatible servers.
We must then show that the coupling is preserved under every possible transition. 
Although this approach can be made to work for particular small network topologies, such as the $N$-model, it is difficult to generalize.
This is due to the complex dynamics of the waiting time threshold policy: we find that customers who are routed to a certain server in the original system may be routed to a different server in the upper bound system. 
Thus, the pathwise coupling fails at the level of individual customer types. 
One could proceed by showing a coupling on the level of servers, for example, by proving that whenever a server is busy in the upper bound system, it must be busy in the original system as well, indicating that the original system processes work more efficiently. 
However, we found that even this weaker property fails for most network topologies.

A third method is fluid limit analysis. 
On countable state spaces, fluid limit stability can imply stochastic stability~\cite{Dai1995,Bramson2006,Robert2003}.  
In our model, however, the fluid limit of the number of customers in the system is not informative enough: the departure rate, and hence the slope of the fluid trajectory, depends on whether the waiting times of the customers in front have exceeded their thresholds. 
Slope changes are therefore state-dependent and cannot be inferred from the customer-count trajectory itself.

Instead, we analyze the \gls{FIL} waiting time process to prove stability. 
This process has been used in the analysis of threshold policies, but only in small models with at most two customer types and two servers~\cite{Koole2012, Koole2015, Bekker2011}. 
Our analysis holds for systems with an arbitrary number of customer types and servers. 
We show that the \gls{FIL} waiting time process satisfies the conditions of the \gls{PDMP} framework of \cite{Davis1984}, in particular the boundary condition, which is the distinguishing feature of the PDMP framework compared to standard Markov chain analysis. 
Since the state space is uncountably large, the Foster–Lyapunov criterion \cite{Meyn1993} requires negative drift outside a petite set. 
To this end, we show that bounded waiting time sets are petite by proving that the system can be emptied in finite time with strictly positive probability. 

The main contributions of this paper can thus be summarized as follows. 
First, we establish the exact stability region of skill-based queues with \gls{FCFS} service and waiting-time thresholds. 
Second, we introduce a \gls{PDMP} formulation for the \gls{FIL} waiting time process.
Third, we prove sufficiency of the stability conditions by constructing a Lyapunov function that satisfies the boundary condition of the \gls{PDMP}. 
Since the state space is uncountable, we also prove that bounded waiting-time  sets are petite.

The remainder is organized as follows. 
We introduce the model and \gls{PDMP} formulation in \Cref{sec:model}. 
The stability conditions are proved in \Cref{sec:stab}. 
To illustrate the method, we discuss an illustrative example in \Cref{sec:example}.
Lastly, we discuss the structural assumption of at most one threshold line per server in \Cref{sec:assumption}.

\begin{figure}[h]
    \begin{center}
        \begin{tikzpicture}[
                server/.style={circle, minimum size = .8cm, thick,draw},
                scale=.8
            ]
            \scriptsize
            \foreach \x in {1,2,3}{
                    \node[draw = none] at (1+2.5*\x,0) (queue\x) {};
                    \node[above of = queue\x, yshift = -.5cm] (labelqueue\x) {};
                    \draw[thick] (queue\x.center) --++(-.4,0) --++(0,1.5);
                    \draw[thick] (queue\x.center) --++(.4,0) --++(0,1.5);
                    \node[draw = none, above of = queue\x, yshift = 0cm] {\large $\lambda_\x$};
                }

            \fill[black] ($(queue1.center) + (-.3,.35)$) rectangle ++(.6,-0.2);
            \fill[black] ($(queue2.center) + (-.3,.35)$) rectangle ++(.6,-0.2);
            \fill[black] ($(queue1.center) + (-.3,.65)$) rectangle ++(.6,-0.2);
            \fill[black] ($(queue3.center) + (-.3,.35)$) rectangle ++(.6,-0.2);

            \foreach \x in {1,2,3,4}{
                    \node[server] at (2.4*\x,-2) (server\x) {};
                    \node at (server\x)  {\large $\mu_\x$};
                }

            \draw[very thick, black] (queue1.center) -- (server1.north);
            \draw[very thick, black, dashed] (queue1.center) -- (server2.north) node[pos=0.8, left,yshift=-.05cm]
            {\large $\tau_{12}$};
            \draw[very thick, black, dashed] (queue1.center) -- (server4.north) node[pos=0.6, right, yshift = .15cm]
            {\large $\tau_{14}$};
            \draw[very thick, black] (queue2.center) -- (server1.north);
            \draw[very thick, black] (queue2.center) -- (server2.north);
            \draw[very thick, black, dashed] (queue2.center) -- (server3.north) node[pos=0.8, left,yshift=-.05cm]
            {\large$\tau_{23}$};
            \draw[very thick, black] (queue3.center) -- (server4.north);
        \end{tikzpicture}
    \end{center}
    \caption{A skill-based queue with waiting time thresholds on compatibility lines $(12)$, $(23)$, and $(14)$. The waiting type-2 customer can only start service at server~3 if its waiting time exceeds $\tau_{23}$. }
    \label{fig:intro_qsystem_example}
\end{figure}

\section{Model and PDMP formulation}
\label{sec:model}
Let $\calG = (\calI\cup \calJ, \calL)$ be a bipartite graph with customer types $\calI$, servers $\calJ$, and compatibility lines $\calL \subseteq \calI\times\calJ$.
Type-$i$ customers arrive at queue~$i$ according to an independent Poisson process with rate $\lambda_i \in \R_{>0}$, and service times at server~$j$ are independent and exponentially distributed with rate~$\mu_j\in\R_{>0}$. 
Type-$i$ customers can only be served by server~$j$ if $(ij)\in\calL$.
Let $\calL_T\subseteq \calL$ denote the set of threshold lines, and let $\calL_T^c := \calL \setminus \calL_T$ denote its complement. 
For $(ij)\in\calL_T$, type-$i$ customers are eligible for server~$j$ if and only if their waiting time exceeds $\tau_{ij} \in \R_{> 0}$.  

Service is \gls{FCFS} over all eligible matches, i.e., if a server becomes available then we assign the customer with the largest waiting time among all eligible customers. 
If a customer arrives while there are two or more eligible idle servers, we break ties by assigning the server with the smallest index. 
Similarly, if two or more threshold lines of the same customer type become active at the same time, we break ties by choosing the server with the smallest index.
These tie-breaking rules are used for notational convenience: our analysis can be extended to allow for random tie-breaking as well. 

We write bold letters for vectors and we write $\bfx^{i\mapsto v}$ for the vector $\bfx$ with the $i$-th element replaced by $v$.

\subsection{FIL waiting time process}
We study the waiting time of the \gls{FIL} customers in the system. 
For $i\in\calI$, let $W_i(t) \in \R_{\geq 0}$ denote the waiting time of the \gls{FIL} customer in queue~$i$ at time $t$, where $W_i(t)=0$ if queue~$i$ is empty. 
Suppose that $W_i(t_0)=w_i$ and that the \gls{FIL} type-$i$ customer starts service at time $t_1>t_0$, then $W_i$ increases linearly on $(t_0,t_1)$, i.e.,
$W_i(s)=w_i+s-t_0$ for $s\in(t_0,t_1)$.
At time $t_1$, the next type-$i$ customer, if any, becomes \gls{FIL}. 
Since interarrival times are exponential, the waiting time jumps down to
$
    W_i(t_1)=(W_i(t_1^-)-A_i)^+,
$
where $A_i\sim\expdist{\lambda_i}$ and $x^+:=\max\{x,0\}$.
See \Cref{fig:W_process_N_model} for an illustration of this process for the $N$-model with a threshold on the diagonal line.

Note that the process does not record queue lengths. 
Arrivals to non-empty queues are invisible to the process: such arrivals join behind the current \gls{FIL} customer and do not affect $W_i(t)$ until they themselves become \gls{FIL}.

\subsection{Piecewise deterministic Markov process}
Next, following \cite{Davis1984}, we formulate a \gls{PDMP} $\bfX(t) = (\bfE,\bfW)(t)$, where $\bfE(t)$ models the mode of the system and $\bfW(t)$ is the waiting time vector of all \gls{FIL} customers at time $t \geq 0$. 
The mode records which queues are non-empty and which servers are busy. 
Let
\begin{align}
    \calK &:= \bigl\{ (\bfq,\bfs): 
    \ 
    q_i \in \{+,-\},
    \
    s_j \in\{0,1\}, 
    \ 
    i \in\calI, 
    \ 
    j\in\calJ\bigr\}
    .
\end{align}
Here, $q_i=+$ indicates that queue $i$ is non-empty and $q_i=-$ that it is empty, while $s_j=1$ indicates that server $j$ is busy and $s_j=0$ that it is idle.

We specify the range for $\bfw$ for each mode in $ \calK$. 
If queue~$i$ is empty, then necessarily $w_i=0$.
If queue~$i$ is non-empty and there is an idle server $j\in\calJ$ with threshold $\tau_{ij}$, then the waiting time of the \gls{FIL} type-$i$ customer cannot exceed $\tau_{ij}$.
For $\bfe = (\bfq,\bfs)\in\calK$ and $i\in\calI$ define
\begin{align}
    M_i(\bfe)
    :=
    \begin{cases}
        (0,\bar{\tau}_i(\bfe)) &\textrm{if} \ q_i=+, \\
        \{0\}, &\textrm{if} \ q_i=-,
    \end{cases}
    \qquad 
    \textrm{where} 
    \
    \bar{\tau}_i(\bfe) 
    := 
    \inf\{\tau_{ij} : (ij)\in\calL_T,\ s_j=0\}
    ,
\end{align}
with the convention that the infimum over the empty set equals $+\infty$. 
Here, $\bar{\tau}_i(\bfe)$ is the smallest threshold from type~$i$ to an idle server.
The state space of $\bfX(t)$ is 
\begin{align}
    E
    &:=
    \{(\bfe,\bfw): \bfe\in\calK,\ \bfw\in M_{\bfe}\}
    ,
    \qquad 
    \textrm{with} 
    \ 
    M_{\bfe}
    :=
    \prod_{i\in\calI} M_i(\bfe)
    . 
\end{align}

The probability law of $\bfX(t)$ is determined by three ingredients:
(i) a deterministic flow $\phi_\bfe(t,\bfw)$, governing the motion in mode $\bfe\in\calK$ between jumps,
(ii) an interior jump rate $\Lambda: E \to \R_{\geq 0}$, and
(iii) a jump kernel $Q: (E\cup\Gamma) \times \scrB(E) \to [0,1]$, where $\scrB(E)$ is the Borel $\sigma$-field on $E$, and $\Gamma$ is the active boundary of $E$, which we define next. 
Let 
$
    \Gamma
    :=
    \bigcup_{(ij)\in\calL_T} \Gamma_{ij}
$, with
\begin{align}
    \label{eq:def_bnd}
    \Gamma_{ij}
    :=
    \bigl\{
    (\bfq,\bfs,\bfw)\in E:
    q_i=+,\ 
    w_i=\bar{\tau}_i(\bfq,\bfs),\
    j = \min(k: \tau_{ik} = \bar{\tau}_i(\bfq,\bfs))
    \bigr\}
    , 
    \ \
    (ij)\in\calL_T
    . 
\end{align}
In other words, $\Gamma_{ij}$ is the set of states where the \gls{FIL} waiting time of customer type $i$ reaches a threshold value, and server~$j$ is selected by the boundary tie-breaking rule. 
$\Gamma$ is the collection over all such states. 

Between jumps, the state evolves according to $\phi_\bfe$ until either an \emph{interior jump} occurs, triggered by arrivals or departures, or the trajectory hits the boundary $\Gamma$, at which point a \emph{boundary jump} occurs. 
In both cases, the post-jump state is drawn from $Q(\bfx,\cdot)$, where $\bfx \in E$ is the pre-jump state.

We next describe $\phi$, $\Lambda$, $Q$, and the dynamics of $\bfX(t)$ in more detail.

\paragraph{Deterministic flow}
The waiting time of the \gls{FIL} customer of non-empty queues increases linearly between jumps. 
Hence, for $\bfe = (\bfq,\bfs)$ the drift function $\phi_\bfe(t,\cdot): M_\bfe \to M_\bfe$ is given by
\begin{align}
    \label{eq:determ_motion_general}
    \phi_{\bfe}(t,\bfw) 
    &= 
    \bigl( \ind{q_i = +}(w_i+t) \bigr)_{i\in\calI}
    .
\end{align}

\paragraph{Interior jump rate}
Interior jumps are triggered by arrivals and departures at busy servers. 
Arrivals to non-empty queues do not change the \gls{FIL} waiting time state. 
Hence, the total jump rate in mode $\bfe = (\bfq,\bfs) \in \calK$ is
\begin{align}
    \label{eq:def_rate_function}
    \Lambda(\bfe)
    :=
    \sum_{i\in\calI: q_i=-} \lambda_i
    +
    \sum_{j\in\calJ: s_j=1} \mu_j 
    .
\end{align}

\paragraph{Jump kernel}
On the interior, the jump kernel is given by
\begin{align}
    \label{eq:def_Q_interior}
    Q(\bfx,B)
    &= \! \!
    \sum_{i\in\calI: q_i=-}
    \frac{\lambda_i}{\Lambda(\bfe)}
    Q_i^{\mathrm{arr}}(\bfx,B)
    + \! \!
    \sum_{j\in\calJ: s_j=1}
    \frac{\mu_j}{\Lambda(\bfe)}
    Q_j^{\mathrm{dep}}(\bfx,B),
    \ \  
    B \in \scrB(E), \ 
    \bfx=(\bfe,\bfw)\in E 
    .
\end{align}
Here, $Q_i^{\mathrm{arr}}(\bfx, \,\cdot\,)$ is the kernel corresponding to a type-$i$ customer arrival, while $Q_j^{\mathrm{dep}}(\bfx, \,\cdot\,)$ is the kernel corresponding to a departure at server~$j$, which we specify next.

{\em Arrivals.}
Suppose that the system is in state $\bfx=(\bfq,\bfs,\bfw)\in E$ with $q_i=-$ when a type-$i$ customer arrives. 
Let
$
    \calS_i(\bfx)
    :=
    \{j\in\calJ: s_j=0,\ (ij)\in\calL_T^c \}
$
denote the set of idle servers that can serve an arriving type-$i$ customer immediately.
If \mbox{$\calS_i(\bfx)\neq\emptyset$}, then the arriving customer is routed to server
$\ell_i(\bfx) := \min \calS_i(\bfx)$. 
On the other hand, if $\calS_i(\bfx)=\emptyset$ then the customer cannot  be served immediately, so it joins queue~$i$.
We thus have
\begin{align}
    \label{eq:def_Q_int_arrival}
    Q_i^{\mathrm{arr}}(\bfx,B)
    &=
    \begin{cases}
        \ind{(\bfq^{i \mapsto +},\bfs,\bfw)\in B},
        & \calS_i(\bfx)=\emptyset
        ,
        \\
        \ind{(\bfq,\bfs^{\ell_i(\bfx) \mapsto 1},\bfw)\in B},
        & \calS_i(\bfx)\neq\emptyset
        ,
    \end{cases}
    \qquad 
    B\in\scrB(E)
    .
\end{align}

{\em Departures.}
Suppose that the system is in state $\bfx=(\bfq,\bfs,\bfw)\in E$ with $s_j=1$ when there is a departure at server~$j$. 
Let
$
    \calC_j(\bfx)
    :=
    \{i\in\calI: q_i=+,\ (ij)\in\calL,\ w_i\geq \tau_{ij} \ind{(ij)\in\calL_T}\}
$
denote the set of \gls{FIL} customers that are eligible to start service at server~$j$.
If $\calC_j(\bfx) \neq\emptyset$, then, by the \gls{FCFS} policy, the \gls{FIL} customer of type 
$
    m_j(\bfx)
    :=
    \argmax_{i\in\calC_j(\bfx)} w_i
$
starts service.
If $\calC_j(\bfx)=\emptyset$, then no eligible customer is available, hence server~$j$ becomes idle.

Suppose that the \gls{FIL} type~$m$ customer with waiting time $w_m$ starts service at server~$j$. 
Let $A_m$ denote the interarrival time between this customer and its successor, with cumulative distribution function $F_{A_m}(x) := 1 - \exp(-\lambda_m x)$.
If $A_m>w_m$, then no successor has arrived yet and queue~$m$ becomes empty.
If $A_m\leq w_m$, then the successor is already waiting and becomes the new \gls{FIL} customer, whose waiting time is $w_m-A_m$.
We thus have
\begin{align}
    \label{eq:def_Q_int_departure}
    Q_j^{\mathrm{dep}}(\bfx,B)
    &=
    \begin{cases}
        \ind{(\bfq,\bfs^{j\mapsto 0},\bfw)\in B},
        & \calC_j(\bfx)=\emptyset 
        ,
        \\
        c_{m_j(\bfx),j}(\bfx,B),
        & \calC_j(\bfx)\neq\emptyset
        ,
    \end{cases}
    \qquad 
    B\in\scrB(E)
    .
\end{align}
Here,  $c_{m,j}(\bfx)$ is the post-jump distribution when a type-$m$ customer starts service at server~$j$: 
\begin{align}
    \label{eq:def_c_line}
    c_{m,j}(\bfx,B)
    &=
    \bigl(1-F_{A_m}(w_m)\bigr)
    \ind{
        \bigl(
            \bfq^{m\mapsto -},
            \bfs,
            \bfw^{m\mapsto 0}
        \bigr)
        \in B
    }
    \nonumber \\
    &\quad 
    +
    \int_0^{w_m}
    \ind{
        \bigl(
            \bfq,
            \bfs,
            \bfw^{m\mapsto w_m-a}
        \bigr)
        \in B
    }
    \rmd F_{A_m}(a)
    .
\end{align}

{\em Boundary kernel.} Finally, we need to specify the boundary kernel. 
If the process hits $\bfx=(\bfq,\bfs,\bfw)\in \Gamma_{ij}$, then the \gls{FIL} type-$i$ customer starts service at server~$j$. 
Since $w_i=\tau_{ij}$ on $\Gamma_{ij}$, the boundary kernel is given by
\begin{align}
    \label{eq:def_Q_bnd}
    Q(\bfx, B)
    &=
    \bigl(1-F_{A_i}(\tau_{ij})\bigr)\,
    \ind{
        \bigl(
            \bfq^{i\mapsto -},
            \bfs^{j\mapsto 1},
            \bfw^{i\mapsto 0}
        \bigr)\in B
    }
    \nonumber \\
    &\quad
    +
    \int_0^{\tau_{ij}}
    \ind{
        \bigl(
            \bfq,
            \bfs^{j\mapsto 1},
            \bfw^{i\mapsto \tau_{ij}-a}
        \bigr)\in B
    }
    \rmd F_{A_i}(a),
    \qquad
    (\bfq,\bfs,\bfw)\in \Gamma_{ij}.
\end{align}

\paragraph{Process dynamics}
Starting from $\bfX(0)=(\bfe,\bfw)\in E$, the process follows the deterministic flow $\phi_\bfe(\cdot,\bfw)$ until the first jump time $T_1$, at which point it jumps to a new state drawn from $Q(\bfX(T_1^{-}), \, \cdot \,)$.
Jump time $T_1$ is the minimum of the next interior jump time and the next boundary jump time, i.e., 
\begin{align}
\label{eq:survivor_function}
    \P(T_1 > t \mid \bfX(0) = (\bfe,\bfw)) = 
    \begin{cases}
        \exp\bigl( - \Lambda(\bfe)t \bigr), 
        & t < t^*(\bfe,\bfw), \\
        0, 
        & t \geq t^*(\bfe,\bfw)
        .
    \end{cases}
\end{align}
Here,
\begin{align}
\label{eq:t_star_def}
    t^*(\bfe,\bfw) 
    &:= \inf\bigl\{ \tau_{ij} - w_i \,:\, (ij)\in\calL_T,\ q_i = +,\ s_j = 0,\ w_i < \tau_{ij}\bigr\},
    \qquad 
    \bfe = (\bfq,\bfs)
\end{align}
is the smallest time at which the flow hits the active boundary $\Gamma$, which is infinite if all waiting times exceed their thresholds. 

The trajectory of $\bfX(t)$ for $t\leq T_1$ is given by
\begin{align}
    \label{eq:def_dynamics_PDMP}
    \bfX(t) &= 
    \begin{cases}
        \bigl(\bfe, \phi_\bfe(t,\bfw)\bigr), & t < T_1, \\
        (\bfE_1, \bfZ_1), & t = T_1,
    \end{cases}
\end{align}
where $(\bfE_1, \bfZ_1) \sim Q\bigl((\bfe, \phi_\bfe(T_1,\bfw)), \, \cdot \,\bigr)$.
The next jump times $T_2,T_3,\dots$ and post-jump locations $\bfX(T_2),\bfX(T_3),\dots$ are constructed similarly.
The author in \cite{Davis1984} shows that $\bfX(t)$ in~\eqref{eq:def_dynamics_PDMP} is a strong Markov process. 
Let $\scrA_\bfX$ denote the extended generator of $\bfX$ with domain $\calD(\scrA_\bfX)$.

\section{Stability conditions}
\label{sec:stab}
We show in \Cref{thm:stab} that the necessary and sufficient stability conditions for the threshold policy are given by 
\begin{align}
    \label{eq:stab_cond_general}
    \sum_{i\in \calA}\lambda_i
    <
    \sum_{j\in \calS(\calA)}\mu_j,
    \qquad
    \forall \, \emptyset \neq \calA\subseteq \calI
    ,
\end{align}
where $\calS(\calA):=\{j\in\calJ:\exists i\in \calA \text{ with } (ij)\in\calL\}$ is the set of servers that can serve customers of types in~$\calA$.
The authors in~\cite{Adan2014} show that \eqref{eq:stab_cond_general} also characterize the stability region for skill-based queues under \gls{FCFS} service without thresholds.
Our main result establishes that waiting time thresholds do not affect the stability region. 

Our analysis does restrict to skill-based queues where each server has at most one threshold line: 

\begin{assumption}
\label{ass:one_positive_threshold_per_server}
We have
$
    |\{\, i\in\calI : (ij)\in\calL_T\,\}| \leq 1
$
for all
$j\in\calJ$.
\end{assumption}

\begin{theorem}\label{thm:stab} 
    Suppose that \Cref{ass:one_positive_threshold_per_server} holds.
    Then $\bfX(t)$ is positive recurrent if and only if~\eqref{eq:stab_cond_general} holds.
\end{theorem}

We outline the proof of \Cref{thm:stab}; the full proof is in \Cref{app:proof_thm_stab}.
Sufficiency is proved using the Foster--Lyapunov criterion~\cite{Meyn2009}, which we do in three steps.
Firstly, we construct a Lyapunov function $f$ in the domain of the extended generator and evaluate its drift explicitly (\Cref{lem:generator}).
Next, we show that the drift of $f$ is proportional to a weighted sum of waiting times (\Cref{lem:general_drift}).
Lastly, we verify that bounded waiting time sets are petite (\Cref{lem:petite}).
Together, \Cref{lem:general_drift} and \Cref{lem:petite} imply that the drift of $f$ is negative outside a petite set. 
Positive recurrence then follows from the Foster--Lyapunov criterion. 

\Cref{lem:generator} relies on the characterization of the extended generator of \glspl{PDMP} in \cite[Thm.~5.5]{Davis1984}. 
In particular, the author shows that $f$ is in the domain of the extended generator if (i) it is absolutely continuous along trajectories, (ii) the jump differences are integrable, and (iii) it satisfies the boundary condition. 
The boundary condition is characteristic for \glspl{PDMP}.
It states that at a boundary point $\bfx\in\Gamma$, the value $f(\bfx)$ must equal the expected post-jump value $\int_E f(\bfy)\,Q(\bfx,\rmd\bfy)$. 
Our Lyapunov function is constructed precisely to satisfy this boundary condition.

We next introduce the Lyapunov function and formalize the proof steps.
Let 
\begin{align}
    \label{eq:def_idle_server_set}
    \calJ_i^T(\bfs) 
    &:= 
    \{j\in\calJ: (ij)\in\calL_T, \ s_j=0 \},
    \qquad 
    i\in\calI, \ 
    \bfs \in \{0,1\}^{|\calJ|}
\end{align}
be the set of idle servers connected to $i\in\calI$ via threshold lines.
Moreover, let
\begin{align}
    \label{eq:def_kappa_ij_star}
    \kappa_{ij}^*
    &:=
    \frac{1}{\tau_{ij}^2}
    \int_0^{\tau_{ij}}(\tau_{ij}-a)^2\rmd F_{A_i}(a)
    \in (0,1)
    ,
    \qquad 
    (ij)\in\calL_T
\end{align}
be the normalized expected squared residual waiting time that remains after a boundary jump.

We consider the Lyapunov function
\begin{align}
    \label{eq:def_test_function}
    f(\bfq,\bfs,\bfw) 
    &:= 
    \sum_{i\in\calI} \alpha_i(\bfs) w_i^2,
    \qquad
    \textrm{where} 
    \
    \alpha_i(\bfs) 
    := 
    \lambda_i
    \prod_{j\in\calJ_i^T(\bfs)} \kappa^*_{ij},
    \qquad 
    (\bfq,\bfs,\bfw) \in E
    .
\end{align}
Here, we use the convention that the product over an empty set equals one.

The coefficients $\alpha$ in~\eqref{eq:def_test_function} are constructed so that 
\begin{align}
    \label{eq:alpha_property}
    \alpha_i(\bfs)
    &= 
    \begin{cases}
        \kappa_{ij}^*  \alpha_i(\bfs^{j\mapsto 1}), &\textrm{if} \ j\in \calJ_i^T(\bfs), \\
        \alpha_i(\bfs^{j\mapsto 1}), &\textrm{else}
        .
    \end{cases} 
\end{align}
Indeed, if $j\in \calJ_i^T(\bfs)$ then $\calJ_i^T(\bfs^{j\mapsto 1}) = \calJ_i^T(\bfs)\setminus\{j\}$ by~\eqref{eq:def_idle_server_set} so the product in~\eqref{eq:def_test_function} changes by one element. 
If $j\notin \calJ_i^T(\bfs)$ then $\calJ_i^T(\bfs^{j\mapsto 1}) = \calJ_i^T(\bfs)$ and so $\alpha_i(\bfs) = \alpha_i(\bfs^{j\mapsto 1})$.
Under \Cref{ass:one_positive_threshold_per_server}, for each fixed server $j$ there is at most one type~$i$ such that $j\in\calJ_i^T(\bfs)$. 
Hence, changing server $j$ from idle to busy can affect at most one coefficient $\alpha_i$ in the test function.

We evaluate the drift of $f$ in \Cref{lem:generator}.
Its proof in \Cref{app:proof_lem_generator} follows by verifying that $f$ satisfies the conditions of \cite[Thm.~5.5]{Davis1984}.

\begin{lemma}
\label{lem:generator}
    If \Cref{ass:one_positive_threshold_per_server} holds,
    then $f \in \calD(\scrA_\bfX)$ and 
    \begin{align}
        \label{eq:drift_decomposition}
        \scrA_\bfX f(\bfx)
        =
        \sum_{i\in\calI} 2\alpha_i(\bfs) w_i
        +
        \Lambda(\bfe) \int_E\bigl(f(\bfy)-f(\bfx)\bigr) Q(\bfx, \rmd\bfy)
        ,
        \qquad 
        \bfx=(\bfe,\bfw)\in E, \
        \bfe=(\bfq,\bfs) \in \calK
        .
    \end{align}
\end{lemma}

The second step in the sufficiency proof is to show that $f$ has drift proportional to the weighted waiting times in  \Cref{lem:general_drift}. 
Its proof is in \Cref{app:proof_lem_general_drift}.

\begin{lemma}\label{lem:general_drift}
    If \Cref{ass:one_positive_threshold_per_server} and~\eqref{eq:stab_cond_general} hold, then there exist $\eta,c_0 >0$ such that
    \begin{align}
        \scrA_\bfX f(\bfe,\bfw) 
        &\leq 
        -\eta
        \sum_{i\in\calI} \alpha_i(\bfs) w_i + c_0,
        \qquad 
        (\bfe,\bfw) \in E, \ \bfe = (\bfq,\bfs) \in \calK
        .
    \end{align}
\end{lemma}

Recall that a measurable set $\calE \subset E$ is \emph{petite} if there exist a probability measure $\gamma$ on $\R_{\geq 0}$ and a non-trivial measure $\nu$ on $\scrB(E)$ such that for all $\bfx\in\calE$ and $B\in\scrB(E)$,
$
    \int_0^\infty \P_\bfx(\bfX(t)\in B)\,\gamma(\rmd t)
    \geq
    \nu(B)
$
\cite{Meyn2009}.
Moreover, $\calE$ is \emph{small} if there exist $t>0$ and a non-trivial measure $\nu$ on $\scrB(E)$ such that for all $\bfx \in \calE$ and $B\in\scrB(E)$, $\P_\bfx(\bfX(t)\in B) \geq \nu(B)$. 
As shown in~\mbox{\cite[Sec.~5.5.2]{Meyn2009}}, small sets are petite.
\Cref{lem:petite} shows that bounded waiting time sets are small, and therefore petite. 
The proof in \Cref{app:proof_lem_petite} shows that when starting from a state with bounded waiting times, there is a positive probability that the system becomes empty in finite time. 
This proves that bounded waiting time sets are small with respect to the Dirac measure at the empty state.

\begin{lemma}\label{lem:petite}
    For any $a \geq 0$,
    $
        B_a
        :=
        \{(\bfe,\bfw)\in E: w_i\leq a, \ \forall \, i \in \calI\}
    $
    is small.
\end{lemma}

The sufficiency proof in \Cref{app:proof_thm_stab} is concluded from \Cref{lem:general_drift}, \Cref{lem:petite}, and the Foster--Lyapunov criterion.

Necessity of \eqref{eq:stab_cond_general} follows from a standard coupling argument: suppose that there exists $\calA\subseteq \calI$ for which \eqref{eq:stab_cond_general} fails. 
We construct a coupled single server queue whose queue length is always at most the total number of customers of types in $\calA$ in the original system.
We then show that the single server queue is unstable, and therefore the original process $\bfX(t)$ cannot be positive recurrent.
See \Cref{app:proof_thm_stab} for the full proof. 

\section{Illustrative example: the $N$-model}
\label{sec:example}
We illustrate the \gls{PDMP} construction and the Lyapunov function for the $N$-model \cite{Gans2003}, see \Cref{fig:W_process_N_model}. 
Let $\calI = \calJ = \{1,2\}$, $\calL = \{(11),(12),(22)\}$ and $\calL_T = \{(12)\}$ with threshold $\tau_{12} \in \R_{> 0}$. 
\Cref{fig:W_process_N_model} also shows a possible realization of the \gls{FIL} waiting time process.

For this system, the set of modes is $\calK = \{(q_1,q_2,s_1,s_2), \ q_i \in\{+,-\}, \ s_j \in \{0,1\}, \ \forall \, i,j\in \{1,2\}\}$. 
In mode $\bfe^* := (+,-,1,0)$, queue~1 is non-empty and server~2 is idle.
For this mode, $W_1(t) <\tau_{12}$. 
Indeed, if $W_1(t) \geq \tau_{12}$ then the \gls{FIL} type-1 customer is eligible for server~2, and so server~2 cannot be idle. 
The only active boundary is thus $\Gamma = \Gamma_{12} = \{(\bfe^*,\tau_{12},0)\}$.
The Lyapunov function in~\eqref{eq:def_test_function} evaluates to 
\begin{align}
    \label{eq:test_function_N}
    f(\bfe,\bfw)
    &= 
    \begin{cases}
        \lambda_1 \kappa_{12}^* w_1^2 + \lambda_2 w_2^2, &\textrm{if} \ \bfe = \bfe^*, \\  
        \lambda_1w_1^2 + \lambda_2 w_2^2,  &\textrm{else}
        ,
    \end{cases}
\end{align}
with $\kappa_{12}^*$ as in~\eqref{eq:def_kappa_ij_star}.

It is straightforward to verify that $f$ in~\eqref{eq:test_function_N} satisfies the absolute continuity and integrability conditions in \cite[Thm.~5.5]{Davis1984}. 
We next check that it also satisfies the boundary condition.

At the boundary, the \gls{FIL} type-1 customer reaches waiting time $\tau_{12}$ and therefore becomes eligible for server~2. 
The boundary jump corresponds to assigning this customer to server~2. 
If no other type-1 customer has arrived during the waiting time interval of length $\tau_{12}$, then queue~1 becomes empty after the jump and so the post-jump state is $(-,-,1,1,0,0)$. 
Otherwise, the next type-1 customer becomes the new \gls{FIL} customer.
If this customer arrived $a$ time units after the previous \gls{FIL} customer, its waiting time after the jump is $\tau_{12}-a$, and so the post-jump state is $(+,-,1,1,\tau_{12}-a)$. 

Letting $\bfz^* := (\bfe^*, \tau_{12},0)$, we thus obtain from~\eqref{eq:def_Q_bnd},
\begin{align}
    \int_E f(\bfy) \, Q(\bfz^*,\rmd\bfy)
    &=
    (1-F_{A_1}(\tau_{12}))
    f(-,-,1,1,0,0)
    +
    \int_0^{\tau_{12}} \! \!
    f(+,-,1,1,\tau_{12}-a,0)
    \rmd F_{A_1}(a) \\
    &\stackrel{\eqref{eq:test_function_N}}{=}
    0 + \lambda_1 \int_0^{\tau_{12}} (\tau_{12} - a)^2  \rmd F_{A_1}(a) \\
    &\stackrel{\eqref{eq:def_kappa_ij_star}}{=} \lambda_1 \kappa_{12}^* \tau_{12}^2 
    \stackrel{\eqref{eq:test_function_N}}{=}
    f(\bfz^*)
    .
    \label{eq:bnd_N_model}
\end{align}

That is, the definition of $\kappa_{12}^*$ is precisely chosen so that~\eqref{eq:bnd_N_model} holds.

\begin{figure}[H]
    \centering
    \begin{tikzpicture}[
                server/.style={circle, minimum size = .7cm, thick,draw},
                scale=.9
            ]

        \foreach \x in {1,2}{
            \node[draw = none] at (-4.2+1*\x,2) (queue\x) {};
            \node[above of = queue\x, yshift = -.5cm] (labelqueue\x) {};
            \draw[thick] (queue\x.center) --++(-.3,0) --++(0,1.5);
            \draw[thick] (queue\x.center) --++(.3,0) --++(0,1.5);
            \node[draw = none, above of = queue\x, yshift = 0cm] { $\lambda_\x$};
        }

        \fill[black] ($(queue1.center) + (-.25,.35)$) rectangle ++(.5,-0.2);
        \fill[black] ($(queue1.center) + (-.25,.65)$) rectangle ++(.5,-0.2);
        \fill[black] ($(queue2.center) + (-.25,.35)$) rectangle ++(.5,-0.2);

        \foreach \x in {1,2}{
            \node[server] at (-4.2+1*\x,0) (server\x) {};
            \node at (server\x)  { $\mu_\x$};
        }

        \draw[very thick, black] (queue1.center) -- (server1.north);
        \draw[very thick, black, dashed] (queue1.center) -- (server2.north) node[pos=0.3, right, xshift=-.05cm]
        {\small $\tau_{12}$};
        \draw[very thick, black] (queue2.center) -- (server2.north);

        \draw[->, thick] (-0.5,0) -- (15,0) node[below] {$t$};
        \draw[->, thick] (0,-0.5) -- (0,5) node[left] {$W_i(t)$};
        \draw[thick,dashed] (0,3.5) node[left] {$\tau_{12}$}  -- (15,3.5); 

        \draw[thick,blue] (2,0) -- (4,2) node[pos=0.5, above, rotate=45] {\small $W_1(t)$};
        \draw[thick,dashed,blue] (4,2) -- (4,.5);
        
        \draw[thick,blue] (4,.5) -- (7,3.5);
        \draw[thick,dashed,blue] (7,3.5) -- (7,1); 
        \draw[thick,blue] (7,1) -- (9,3);
        \draw[thick,dashed,blue] (9,3) -- (9,2.25); 
        \draw[thick,blue] (9,2.25) -- (11.25,4.5);
        \draw[thick,dashed,blue] (11.25,4.5) -- (11.25,.5); 
        \draw[thick,blue] (11.25,.5) -- (13,2.25);
        \draw[thick,dashed,blue] (13,2.25) -- (13,0); 

        \draw[ultra thick,red] (8,0) -- (9.7,1.7) node[pos=0.5, above, rotate=45] {\small $W_2(t)$};;
        \draw[ultra thick,dashed,red] (9.7,1.7) -- (9.7,.9);
        \draw[ultra thick,red] (9.7,.9) -- (11.9,3);
        \draw[ultra thick,dashed,red] (11.9,3) -- (11.9,0);

        \def\yA{-0.7}   
        \def\yB{-1.2}   
        \def\yC{-1.7}   
        \def\yD{-2.2}   

        \draw[gray!50] (0,\yA) -- (15,\yA);
        \draw[gray!50] (0,\yB) -- (15,\yB);
        \draw[gray!50] (0,\yC) -- (15,\yC);
        \draw[gray!50] (0,\yD) -- (15,\yD);

        \node[left] at (0,\yA) {\scriptsize arrivals 1};
        \node[left] at (0,\yB) {\scriptsize arrivals 2};
        \node[left] at (0,\yC) {\scriptsize server 1};
        \node[left] at (0,\yD) {\scriptsize server 2};

        \draw[densely dotted, blue] (1.25,\yA) -- (1.25,0); 
        \node[draw=blue,circle,inner sep=1.2pt,fill=white, thick] at (1.25,\yA) {\scriptsize \color{blue} 1}; 
        \draw[densely dotted, blue] (2,\yA) -- (2,0);  
        \node[draw=blue,circle,inner sep=1.2pt,fill=white, thick] at (2,\yA) {\scriptsize \color{blue}  2}; 
        \draw[densely dotted, blue] (3.5,\yA) -- (3.5,0);  
        \node[draw=blue,circle,inner sep=1.2pt,fill=white, thick] at (3.5,\yA) {\scriptsize \color{blue}  3};
        \draw[
            decorate,
            color=blue,
            decoration={brace,amplitude=3.5pt}
        ] (2.05,\yA+0.22) -- (3.45,\yA+0.22)
        node[midway,above=5pt,fill=white,draw,fill opacity=0.7, text opacity=1,inner sep=1pt,] {\scriptsize $A_{13}$};
        \draw[densely dotted, blue] (6,\yA) -- (6,0); 
        \node[draw=blue,circle,inner sep=1.2pt,fill=white, thick] at (6,\yA) {\scriptsize \color{blue}  4}; 
        \draw[densely dotted, blue] (6.75,\yA) -- (6.75,0);
        \node[draw=blue,circle,inner sep=1.2pt,fill=white, thick] at (6.75,\yA) {\scriptsize \color{blue}  5}; 
        \draw[densely dotted, blue] (10.75,\yA) -- (10.75,0); 
        \node[draw=blue,circle,inner sep=1.2pt,fill=white, thick] at (10.75,\yA) {\scriptsize \color{blue}  6}; 

        \draw[densely dotted, red] (8,\yB) -- (8,0);
        \node[draw=red,circle,inner sep=1.2pt,fill=red, thick] at (8,\yB) {\scriptsize \color{white} \bf 1};
        \draw[densely dotted, red] (8.8,\yB) -- (8.8,0);
        \node[draw=red,circle,inner sep=1.2pt,fill=red, thick] at (8.8,\yB) {\scriptsize \color{white} \bf 2};
        \draw[
            decorate,
            color=red,
            decoration={brace,amplitude=3.5pt}
        ] (8.05,\yB+0.22) -- (8.75,\yB+0.22)
        node[midway,above=5pt,fill=white,draw,fill opacity=0.7, text opacity=1,inner sep=1pt,] {\scriptsize $A_{21}$};

        \draw[very thick,blue] (1.25,\yC) -- (4,\yC)
            node[midway,draw,circle,inner sep=1.2pt,fill=white, thick] {\scriptsize 1};
        \draw[very thick,blue] (4,\yC) -- (9,\yC)
            node[midway,draw,circle,inner sep=1.2pt,fill=white, thick] {\scriptsize 2};
        \draw[very thick,blue] (9,\yC) -- (11.25,\yC)
            node[midway,draw,circle,inner sep=1.2pt,fill=white, thick] {\scriptsize 4};
        \draw[very thick,blue] (11.25,\yC) -- (13,\yC)
            node[midway,draw,circle,inner sep=1.2pt,fill=white, thick] {\scriptsize 5};
        \draw[very thick,blue] (13,\yC) -- (14.7,\yC)
            node[midway,draw,circle,inner sep=1.2pt,fill=white, thick] {\scriptsize 6};

        \draw[very thick,blue] (7,\yD) -- (9.7,\yD)
            node[midway,draw,circle,inner sep=1.2pt,fill=white, thick] {\scriptsize 3};
        \draw[ultra thick,red] (9.7,\yD) -- (11.9,\yD)
            node[midway,draw,circle,inner sep=1.2pt,fill=red, thick] {\scriptsize \color{white} \bf 1};
        \draw[ultra thick,red] (11.9,\yD) -- (12.6,\yD)
            node[midway,draw,circle,inner sep=1.2pt,fill=red, thick] {\scriptsize \color{white} \bf 2};

        \foreach \x in {1.25,4,9,11.25,13,14.7}
            \draw[blue, thick] (\x,\yC+0.14) -- (\x,\yC-0.14);
        \foreach \x in {7,9.7}
            \draw[blue, thick] (\x,\yD+0.14) -- (\x,\yD-0.14);
        \foreach \x in {11.9,12.6}
            \draw[red, thick] (\x,\yD+0.14) -- (\x,\yD-0.14);

        \draw[
            decoration={brace, amplitude=4pt},
            color=blue,
            decorate
            ] (4.07,1.95) --++ (0,-1.3) node[midway, right=5pt, yshift=0cm, fill=white, fill opacity=0.9, text opacity=1, inner sep=1pt, draw] {\scriptsize $A_{13}$};

        \draw[
            decoration={brace, amplitude=3pt},
            color=red,
            decorate
            ] (9.77,1.65) --++ (0,-.6) node[midway, right=4pt, yshift=0cm, fill=white, fill opacity=0.9, text opacity=1, inner sep=1pt, draw] {\scriptsize $A_{21}$};

    \end{tikzpicture}
    \caption{The layout of the $N$-model with threshold on line $(12)$ (left) and a possible realization of the \gls{FIL} waiting time process $(W_1,W_2)(t)$ (right). Arrival times and service durations are shown below the graph, with corresponding customer labels. When $W_1(t)$ hits $\tau_{12}$, the third type-1 customer starts service at server~2. }
    \label{fig:W_process_N_model}
\end{figure}

\section{Possible relaxation of the threshold assumption}
\label{sec:assumption}
We reflect on \Cref{ass:one_positive_threshold_per_server} and the required steps to relax this assumption. 
To be in the domain of the extended generator, the Lyapunov function must satisfy the boundary condition. 
If we let $g(\bfx) = \sum_{i\in\calI} \beta_i(\bfs) w_i^2$ for weights $\beta_i(\bfs)$, $i\in\calI, \bfs\in\{0,1\}^{|\calJ|}$, then the boundary condition for $(ij)\in\calL_T$ evaluates to
\begin{align}
    \label{eq:system_alpha}
    \Bigl( \sum_{m\neq i}\beta_m(\bfs)w_m^2 \Bigr)
    + 
    \beta_i(\bfs)\tau_{ij}^2
    &=
    \Bigl(
    \sum_{m\neq i}
    \beta_m(\bfs^{j\mapsto 1})
    w_m^2
    \Bigr)
    +
    \beta_i(\bfs^{j\mapsto 1})
    \int_0^{\tau_{ij}}(\tau_{ij}-a)^2\rmd F_{A_i}(a)
    .
\end{align}
See \Cref{app:proof_lem_generator} for this derivation.  
Equation \eqref{eq:system_alpha} is satisfied if 
$
    \beta_m(\bfs) = \beta_m(\bfs^{j\mapsto 1})
$
for all $m\neq i$, and 
$
    \beta_i(\bfs) = \kappa_{ij}^* \, \beta_i(\bfs^{j\mapsto 1})
    .
$
Under \Cref{ass:one_positive_threshold_per_server}, these requirements are compatible. 
Indeed,  changing $s_j$ from $0$ to $1$ can affect at most one coefficient.  
This is achieved by the product-form construction in~\eqref{eq:def_test_function}.

Suppose now that there exist $i_1,i_2,j$ such that $(i_1,j),(i_2,j)\in\calL_T$. 
The boundary at $\Gamma_{i_1 j}$ then forces $\beta_{i_2}(\bfs) = \beta_{i_2}(\bfs^{j\mapsto 1})$, while the boundary at $\Gamma_{i_2 j}$ forces $\beta_{i_2}(\bfs) = \kappa_{i_2 j}^* \, \beta_{i_2}(\bfs^{j\mapsto 1})$. 
Since $\kappa_{i_2 j}^* < 1$, these constraints are incompatible.

Relaxing \Cref{ass:one_positive_threshold_per_server} would thus require a structurally different Lyapunov function, as the quadratic form used in our proof does not satisfy the boundary condition for general topologies. Identifying an alternative Lyapunov function that satisfies the boundary condition is known to be technically challenging within the \gls{PDMP} framework, and we therefore leave this direction for future work.

\subsection{Other extensions}
There are several natural extensions of our model. 
For example, we could allow for customer dependent service rates, i.e., type-$i$ customers would be served at server-$j$ with rate $\mu_{ij} > 0$. 
Another extension is to consider priority based routing when the thresholds are exceeded, instead of \gls{FCFS} service.
We believe that the \gls{PDMP} framework itself is sufficiently flexible to accommodate such variants. 
However, the Lyapunov function used in our sufficiency proof is tailored to the \gls{FCFS} dynamics and customer independent service rates.
Extending the drift argument to these settings would therefore require new Lyapunov constructions.

\section*{Acknowledgments}
We are grateful to Jaron Sanders, Ralph van Ierland, Sem Borst, and Maarten Wolf for the valuable discussions and perspectives they offered during the development of this work.

This work is part of \emph{Valuable AI}, a research collaboration between the Eindhoven University of Technology and the Koninklijke KPN N.V.
Parts of this research have been funded by the EAISI's IMPULS program, and by Holland High Tech | TKI HTSM via the PPS allowance scheme for public-private partnerships.

\printbibliography

@article{Meyn1993,
   author = {Sean P Meyn and R L Tweedie},
   doi = {10.2307/1427522},
   issue = {3},
   journal = {Advances in Applied Probability},
   pages = {518-548},
   title = {{Stability of Markovian processes III: Foster-Lyapunov criteria for continuous-time processes}},
   volume = {25},
   year = {1993}
}

@article{Davis1984,
   author = {M H A Davis},
   doi = {10.1111/j.2517-6161.1984.tb01308.x},
   isbn = {46/3/353/7028122},
   issue = {3},
   journal = {Journal of the Royal Statistical Society. Series B (Methodological)},
   pages = {353-388},
   title = {{Piecewise-Deterministic Markov Processes: A general class of non-diffusion stochastic models}},
   volume = {46},
   year = {1984}
}

@article{Adan2009,
   author = {Ivo Adan and Robert D Foley and David R McDonald},
   doi = {10.1007/s11134-009-9140-y},
   issn = {15729443},
   issue = {4},
   journal = {Queueing Systems},
   pages = {311-344},
   title = {{Exact asymptotics for the stationary distribution of a Markov chain: A production model}},
   volume = {62},
   year = {2009}
}

@article{Adan2014,
   author = {Ivo Adan and Gideon Weiss},
   doi = {10.1214/13-ssy117},
   issue = {1},
   journal = {Stochastic Systems},
   pages = {250-299},
   title = {{A skill-based parallel service system under FCFS-ALIS — steady state, overloads, and abandonments}},
   volume = {4},
   year = {2014}
}

@article{Adan2012,
   author = {Ivo Adan and Gideon Weiss},
   doi = {10.1287/opre.1110.1027},
   issn = {0030364X},
   issue = {2},
   journal = {Operations Research},
   pages = {475-489},
   title = {{Exact FCFS matching rates for two infinite multitype sequences}},
   volume = {60},
   year = {2012}
}

@book{Robert2003,
   author = {Philippe Robert},
   doi = {10.1007/978-3-662-13052-0},
   isbn = {9783642056253},
   publisher = {Springer–Verlag},
   title = {Stochastic networks and queues},
   year = {2003}
}

@article{Visschers2012,
   author = {Jeremy Visschers and Ivo Adan and Gideon Weiss},
   doi = {10.1007/s11134-011-9274-6},
   isbn = {1113401192},
   issn = {15729443},
   issue = {3},
   journal = {Queueing Systems},
   pages = {269-298},
   title = {A product form solution to a system with multi-type jobs and multi-type servers},
   volume = {70},
   year = {2012}
}

@article{Dai1995,
   author = {Jim G Dai},
   doi = {10.1214/aoap/1177004828},
   issue = {1},
   journal = {Annals of Applied Probability},
   pages = {49-77},
   title = {{On positive Harris recurrence of multiclass queueing networks: A unified approach via fluid limit models}},
   volume = {5},
   year = {1995}
}

@book{Bramson2006,
   author = {Maury Bramson},
   doi = {10.1007/978-3-540-68896-9},
   isbn = {9783540688952},
   publisher = {Springer-Verlag},
   title = {Stability of queueing networks},
   year = {2006}
}

@article{Koole2015,
   author = {G M Koole and B F Nielsen and T B Nielsen},
   doi = {10.1017/S0269964815000091},
   isbn = {0269964815000},
   issn = {14698951},
   issue = {3},
   journal = {Probability in the Engineering and Informational Sciences},
   pages = {461-471},
   title = {Optimization of overflow policies in call centers},
   volume = {29},
   year = {2015}
}

@article{Bekker2011,
   author = {R Bekker and G M Koole and B F Nielsen and T B Nielsen},
   doi = {10.1007/s11134-011-9225-2},
   isbn = {1113401192},
   issn = {15729443},
   issue = {1},
   journal = {Queueing Systems},
   pages = {61-78},
   title = {Queues with waiting time dependent service},
   volume = {68},
   year = {2011}
}

@article{Koole2012,
   author = {G M Koole and B F Nielsen and T B Nielsen},
   doi = {10.1287/opre.1120.1089},
   issn = {0030364X},
   issue = {5},
   journal = {Operations Research},
   pages = {1258-1266},
   title = {First in line waiting times as a tool for analysing queueing systems},
   volume = {60},
   year = {2012}
}

@book{Meyn2009,
   author = {Sean P Meyn and R L Tweedie},
   doi = {10.1007/978-1-4471-3267-7},
   isbn = {978-0-521-73182-9},
   pmid = {38252772},
   publisher = {Cambridge University Press.},
   title = {Markov chains and stochastic stability},
   year = {2009}
}

@article{Gans2003,
   author = {N Gans and G M Koole and A Mandelbaum},
   doi = {10.1287/msom.5.2.79.16071},
   issn = {15234614},
   issue = {2},
   journal = {Manufacturing and Service Operations Management},
   pages = {79-141},
   title = {Telephone Call Centers: Tutorial, Review, and Research Prospects},
   volume = {5},
   year = {2003}
}

\appendix

\section{Proofs}

\subsection{Proof of \Cref{thm:stab}}
\label{app:proof_thm_stab}

\begin{proof}[Proof of \Cref{thm:stab}]
    {\it Proof of sufficiency.}
    Let $\eps>0$ and let 
    $
        \alpha_{\min}
        :=
        \min_{i, \bfs} \alpha_i(\bfs)
        >0
        ,
    $
    which is well-defined since $\bfs \in \{0,1\}^{|\calJ|}$ and $i \in \calI$, both having finite cardinality. 
    Then by \Cref{lem:general_drift} there exist $\eta,c_0 >0$ such that for all $\bfx\in E$,
    \begin{align}
        \label{eq:general_neg_drift}
        \scrA_\bfX f(\bfx)
        \leq
        -\eps
        +
        (c_0 + \eps)
        \ind{
            \bfx\in\calC
        },
    \end{align}
    with $\calC := \{(\bfe,\bfw)\in E: \ w_i \leq (c_0+\eps)/(\alpha_{\min} \eta), \ \forall \ i\in\calI \}$. 
    Indeed, outside $\calC$ we have $w_i > (c_0+\eps)/(\alpha_{\min}\eta)$ for some $i$, so that $\eta\sum_{i\in\calI} \alpha_i(\bfs)w_i \geq \eta\alpha_{\min} w_i > c_0 + \eps$. 
    It follows from \Cref{lem:petite} that $\calC$ is small and therefore petite~\cite{Meyn2009}.
    Moreover, $f$ is norm-like, i.e., $f(\bfy) \to \infty$ as $\|\bfy\| \to\infty$. 
    Positive recurrence of $\bfX(t)$ follows from \mbox{\cite[Thm.~4.2]{Meyn1993}}.

    {\it Proof of necessity.}
    Suppose that there exists $\calA\subseteq \calI$ such that $\sum_{i\in\calA} \lambda_i \geq \sum_{j\in\calS(\calA)} \mu_j$.
    Let $A(t)$ and $D(t)$ denote the numbers of arrivals and departures, respectively, of customers with types in $\calA$ up to time~$t$.
    Let $P_j(t)$ be a Poisson process with rate $\mu_j$ representing the potential departure process of server~$j\in\calS(\calA)$.
    A jump of $P_j(t)$ results in an actual departure only if server $j$ is busy just before the jump.

    Consider a single server queue with arrival process $A(t)$ and potential departure process $\sum_{j\in\calS(\calA)} P_j(t)$. 
    This is an $M/M/1$ queue with arrival rate $\sum_{i\in\calA} \lambda_i$ and service rate $\sum_{j\in\calS(\calA)} \mu_j$. 
    Let $\bar{D}(t)$ denote the number of departures up to time $t$ in the single server queue. 

    We next show that $D(t) \leq \bar{D}(t)$ for all $t\geq 0$. 
    Suppose that there exists $t_0\geq 0$ such that $D(t_0^-) = \bar{D}(t_0^-)$ and $D(t_0) > \bar{D}(t_0)$, then there is a departure in the skill-based queue but not in the single server queue. 
    Since the potential departure process is the same for both systems, this implies that the single server system is empty at time $t_0$. 
    Hence, all arrivals up to time $t_0$ have already departed from the single server queue and so $\bar{D}(t_0) = A(t_0)$. 
    But then $D(t_0) > \bar{D}(t_0) = A(t_0)$, which is impossible since the number of departures can not exceed the number of arrivals. 
    This shows $D(t) \leq \bar{D}(t)$.

    The inequality $D(t)\leq \bar D(t)$ implies that the skill-based system contains at least
    as many type-$\calA$ customers as the single server queue. 
    Hence, if $\bfX(t)$ were positive recurrent, then the single server queue would also be positive recurrent, contradicting $\sum_{i\in\calA} \lambda_i \geq \sum_{j\in\calS(\calA)} \mu_j$ \cite{Bramson2006}.

\end{proof}

\subsection{Proof of \Cref{lem:generator}}
\label{app:proof_lem_generator}

Before stating the proof of \Cref{lem:generator}, we recall \mbox{\cite[Thm.~5.5]{Davis1984}}, which characterizes the extended generator $\scrA_\bfX$ and its domain $\calD(\scrA_\bfX)$.

\begin{theorem}[Davis]\label{thm:generator}
    Let $g: E\to \R$ be measurable. $g\in \calD(\scrA_\bfX)$ if 
    \begin{enumerate}[label=(P\arabic*), ref=(P\arabic*), noitemsep]
        \item \label{prop:continuous} for any $(\bfe,\bfw)\in E$, the function $t\mapsto g(\bfe,\phi_\bfe(t,\bfw))$ is absolutely continuous for $t\in [0,t^*(\bfe,\bfw))$,
        \item \label{prop:integrability}  for any $t\geq 0$,
        $
            \E\Bigl( \sum_{T_i \leq t} \bigl|g(\bfX(T_i)) - g(\bfX(T_i^-)) \bigr|\Bigr) 
            <
            \infty
            ,
        $
        \item \label{prop:bnd_condition} for any $(\bfe,\bfw)\in \Gamma$, $g$ satisfies the boundary condition
        $
            g(\bfe,\bfw) = \int_{E} g(\bfy) \, Q((\bfe,\bfw), \rmd \bfy)
            .
        $
    \end{enumerate}
    Moreover, $g\in\calD(\scrA_\bfX)$ satisfies
    \begin{align}
        \label{eq:def_extended_gen}
        \scrA_\bfX g(\bfe,\bfw) 
        &= 
        \frac{\rmd}{\rmd t} g(\bfe, \phi_\bfe(t,\bfw)) \Big|_{t=0} + \Lambda(\bfe,w) \int_E \bigl(g(\bfy) - g(\bfe,\bfw)\bigr) \, Q((\bfe,\bfw), \rmd \bfy), \ \ (\bfe,\bfw) \in E
        .
    \end{align}
\end{theorem}

\begin{proof}[Proof of \Cref{lem:generator}]
    We verify that $f$ satisfies \ref{prop:continuous}--\ref{prop:bnd_condition} in \Cref{thm:generator}.

    \ref{prop:continuous}.
    For fixed $\bfe$, $t\mapsto f(\bfe,\phi_\bfe(t,\bfw))$ is a finite sum of quadratic functions, and therefore absolutely continuous. 
    Hence, $f$ satisfies \ref{prop:continuous}.

    \ref{prop:integrability}.
    Let $t\geq 0$ and denote $\lambda_{\max} := \max_{i\in\calI} \lambda_i$.
    Since $\kappa_{ij}^* \in (0,1)$, we have by \eqref{eq:def_test_function} $\alpha_i(\bfs) \leq \lambda_{\max}$ for all $i,j,\bfs$.
    Therefore, $f(\bfe,\bfw) \leq \lambda_{\max} \sum_{i\in\calI} w_i^2$ for all $\bfe, \bfw$. 
    Since $f\geq 0$, it follows that
    \begin{align}
        \E\Bigl(\sum_{T_k \leq t}
        \bigl|f(\bfX(T_k))-f(\bfX(T_k^-))\bigr|\Bigr)
        &\leq
        \E\Bigl(\sum_{T_k \leq t}
        \max\{f(\bfX(T_k)),f(\bfX(T_k^-))\}\Bigr) \\
        &\leq
        \lambda_{\max}
        \E\Bigl(\sum_{T_k \leq t}
        \max\Bigl\{
            \sum_{i\in\calI} W_i(T_k)^2,
            \sum_{i\in\calI} W_i(T_k^-)^2
        \Bigr\}\Bigr).
    \end{align}
    At jumps, the \gls{FIL} waiting times do not increase, while between jumps they
    increase at rate at most one. Hence, for all \(0\leq s\leq t\),
    $
        W_i(s)\leq W_i(0)+s\leq W_i(0)+t 
        .
    $
    Therefore,
    \begin{align}
        \E\Bigl(\sum_{T_k \leq t}
        \bigl|f(\bfX(T_k))-f(\bfX(T_k^-))\bigr|\Bigr)
        &\leq
        \lambda_{\max}
        \E\Bigl(
            \sum_{T_k \leq t}
            \sum_{i\in\calI} (W_i(0)+t)^2
        \Bigr)
        .
    \end{align}
    Let $N(t) = \sum_k \ind{T_k \leq t}$ denote the number of jumps in $[0,t]$,  then
    \begin{align}
        \label{eq:prop_jumps}
        \E\Bigl(\sum_{T_k \leq t}
        \bigl|f(\bfX(T_k))-f(\bfX(T_k^-))\bigr|\Bigr)
        &\leq
        \lambda_{\max}
        \E\bigl(N(t)\bigr)
        \sum_{i\in\calI} (W_i(0)+t)^2
        .
    \end{align}

    Jumps occur due to arrivals, departures, or boundary hits. 
    Since waiting times increase between jumps, there must be a departure between two hits of the same boundary.
    Hence, on every bounded time interval, the number of boundary jumps is bounded by a constant multiple of the number of arrival and service-completion jumps. 
    Since arrivals and departures occur at total rate at most
    $
        \sum_{i\in\calI}\lambda_i+\sum_{j\in\calJ}\mu_j < \infty,
    $
    we have $\E(N(t))<\infty$. 
    Therefore, \eqref{eq:prop_jumps} is finite and  $f$ satisfies \ref{prop:integrability}.

    \ref{prop:bnd_condition}.
    Let $(ij)\in\calL_T$ and let
    $
        \bfx^*=(\bfq,\bfs,\bfw)\in\Gamma_{ij}.
    $
    By definition of $\Gamma_{ij}$ in~\eqref{eq:def_bnd}, $q_i=+$, $s_j=0$, and $w_i=\tau_{ij}$. Substitution into~\eqref{eq:def_test_function} gives
    \begin{align}
        \label{eq:test_x_star_before_quadratic}
        f(\bfx^*)
        &=
        \Bigl( \sum_{m\neq i}\alpha_m(\bfs)w_m^2 \Bigr)
        + 
        \alpha_i(\bfs)\tau_{ij}^2
        .
    \end{align}
    On the other hand, by~\eqref{eq:def_Q_bnd},
    \begin{align}
        \label{eq:boundary_after_quadratic_start}
        \int_E f(\bfy) \, Q(\bfx,\rmd\bfy)
        &=
        \bigl(1-F_{A_i}(\tau_{ij})\bigr)
        \underbrace{f\bigl(
            \bfq^{i\mapsto -},
            \bfs^{j\mapsto 1},
            \bfw^{i\mapsto 0}
        \bigr)}_{:= \mathrm{(i)}}
        +
        \int_0^{\tau_{ij}}
        \underbrace{f\bigl(
            \bfq,
            \bfs^{j\mapsto 1},
            \bfw^{i\mapsto \tau_{ij}-a}
        \bigr)}_{:= \mathrm{(ii)}}
        \rmd F_{A_i}(a).
    \end{align}
    By definition $(\bfw^{i\mapsto 0})_i = 0$ and $(\bfw^{i\mapsto \tau_{ij} - a})_i = \tau_{ij}-a$, hence by~\eqref{eq:def_test_function},
    \begin{align}
        \label{eq:boundary_empty_term_quadratic}
        \mathrm{(i)}
        &=
        \sum_{m\neq i}
        \alpha_m(\bfs^{j\mapsto 1})
        w_m^2 
        , 
        \qquad 
        \mathrm{(ii)}
        =
        \Bigl(\sum_{m\neq i}
        \alpha_m(\bfs^{j\mapsto 1})
        w_m^2 
        \Bigr)
        +
        \alpha_i(\bfs^{j\mapsto 1})
        (\tau_{ij}-a)^2
        .
    \end{align}
    By substituting \eqref{eq:boundary_empty_term_quadratic} into~\eqref{eq:boundary_after_quadratic_start}, the term
    $F_{A_i}(\tau_{ij})
        \sum_{m\neq i}
        \alpha_m(\bfs^{j\mapsto 1})
        w_m^2
    $ 
    cancels and so
    \begin{align}
        \label{eq:boundary_after_quadratic_expanded}
        \int_E f(\bfy) \, Q(\bfx,\rmd\bfy)
        &=
        \Bigl(
        \sum_{m\neq i}
        \alpha_m(\bfs^{j\mapsto 1})
        w_m^2
        \Bigr)
        +
        \alpha_i(\bfs^{j\mapsto 1})
        \int_0^{\tau_{ij}}(\tau_{ij}-a)^2\rmd F_{A_i}(a).
    \end{align}
    It follows from~\eqref{eq:alpha_property} and \Cref{ass:one_positive_threshold_per_server} that
    $
    \alpha_m(\bfs)
        =
        \alpha_m(\bfs^{j\mapsto 1})
    $ for $m\neq i$ and that
    $
    \alpha_i(\bfs)
        =
        \kappa_{ij}^*
        \alpha_i(\bfs^{j\mapsto 1})
    .
    $
    Using this and collecting terms yields
    \begin{align}
        \label{eq:boundary_after_quadratic_simplified}
        \int_E f(\bfy) \, Q(\bfx,\rmd\bfy)
        &=
        \Bigl(
        \sum_{m\neq i}
            \alpha_m(\bfs)
            w_m^2 
        \Bigr)
        +
        \tfrac{1}{\kappa_{ij}^*}
        \alpha_i(\bfs)
        \int_0^{\tau_{ij}}(\tau_{ij}-a)^2\rmd F_{A_i}(a)
        =
        f(\bfx^*)
        .
    \end{align}
    Here, the last equality follows from the definition of $\kappa_{ij}^*$ in~\eqref{eq:def_kappa_ij_star}, and~\eqref{eq:test_x_star_before_quadratic}.

    We have shown that $f$ satisfies \ref{prop:continuous}--\ref{prop:bnd_condition}, therefore $f\in\calD(\scrA_\bfX)$ by \Cref{thm:generator}.\\

    We next evaluate $\scrA_\bfX f$.
    Let $\bfx = (\bfe,\bfw) \in E$ with $\bfe = (\bfq,\bfs)$.
    By~\eqref{eq:def_extended_gen},
    \begin{align}
        \label{eq:extended_gen_exp}
        \scrA_\bfX f(\bfx)
        &=
        \frac{\rmd}{\rmd t} f(\bfe,\phi_{\bfe}(t,\bfw)) \Big|_{t=0}
        +
        \Lambda(\bfe) \int_E\bigl(f(\bfy)-f(\bfx)\bigr) Q(\bfx, \rmd\bfy)
        .
    \end{align}
    By the chain rule, together with~\eqref{eq:determ_motion_general} and~\eqref{eq:def_test_function},
    \begin{align}
        \label{eq:chain_rule}
        \frac{\rmd}{\rmd t} f(\bfe,\phi_{\bfe}(t,\bfw)) \Big|_{t=0}
        =
        \sum_{i\in\calI} 2\alpha_i(\bfs) w_i \ind{q_i=+}
        =
        \sum_{i\in\calI} 2\alpha_i(\bfs) w_i
        ,
    \end{align}
    where the last equality uses that $q_i=-$ implies that $w_i=0$. 
    Substituting~\eqref{eq:chain_rule} into~\eqref{eq:extended_gen_exp} gives~\eqref{eq:drift_decomposition} and completes the proof.
\end{proof}

\subsection{Proof of \Cref{lem:general_drift}}
\label{app:proof_lem_general_drift}

\begin{proof}
    By \Cref{lem:generator} and~\eqref{eq:def_Q_interior},
    \begin{align}
        \label{eq:drift_decomposition_proof}
        \scrA_\bfX f(\bfx)
        =
        \sum_{i\in\calI} 2\alpha_i(\bfs)\, w_i
        +
        \sum_{i\in\calI: q_i=-}
        \lambda_i \Delta_i^{\mathrm{arr}}(\bfx)
        + 
        \sum_{j\in\calJ: s_j=1}
        \mu_j \Delta_j^{\mathrm{dep}}(\bfx)
        ,
    \end{align}
    where
    \begin{align}
        \Delta_i^{\mathrm{arr}}(\bfx)
        &:=
        \int_E
        \bigl(f(\bfy)-f(\bfx)\bigr)
        Q_i^{\mathrm{arr}}(\bfx,\rmd\bfy)
        , 
        \qquad 
        \Delta_j^{\mathrm{dep}}(\bfx)
        :=
        \int_E
        \bigl(f(\bfy)-f(\bfx)\bigr)
        Q_j^{\mathrm{dep}}(\bfx,\rmd\bfy)
        .
    \end{align}

    We first analyze $\Delta_i^{\mathrm{arr}}(\bfx)$. 
    Let $i\in\calI$ with $q_i=-$.
    By~\eqref{eq:def_Q_int_arrival},
    \begin{align}
        \Delta_i^{\mathrm{arr}}(\bfx)
        &=
        \begin{cases}
            f(\bfq^{i\mapsto +},\bfs,\bfw)
            -
            f(\bfq,\bfs,\bfw),
            & \calS_i(\bfx)=\emptyset 
            ,
            \\
            f(\bfq,\bfs^{\ell_i(\bfx)\mapsto 1},\bfw)
            -
            f(\bfq,\bfs,\bfw),
            & \calS_i(\bfx)\neq\emptyset
            .
        \end{cases}
    \end{align}
    Note that
    $f(\bfq^{i\mapsto +},\bfs,\bfw) = f(\bfq,\bfs,\bfw)$, since $\alpha$ in~\eqref{eq:def_test_function} only depends on $\bfs$.
    Therefore, $\Delta_i^{\mathrm{arr}}(\bfx) = 0$ if $\calS_i(\bfx)=\emptyset$. 
    Suppose that $\calS_i(\bfx)\neq\emptyset$, and write
    $
        j=\ell_i(\bfx)
        .
    $
    Then by~\eqref{eq:def_test_function},
    \begin{align}
        \Delta_i^{\mathrm{arr}}(\bfx)
        &=
        \sum_{m\in\calI}
        \bigl(
            \alpha_m(\bfs^{j\mapsto 1})
            -
            \alpha_m(\bfs)
        \bigr)w_m^2 .
    \end{align}
    Since $j\in\calS_i(\bfx)$, we have $s_j=0$. 
    By~\eqref{eq:alpha_property}, $\alpha_m(\bfs) \neq \alpha_m(\bfs^{j\mapsto 1})$ only if $(mj)\in\calL_T$, so 
    \begin{align}
        \Delta_i^{\mathrm{arr}}(\bfx)
        &=
        \sum_{m\in\calI:\,(mj)\in\calL_T}
        \left(
            \tfrac{1}{\kappa_{mj}^*}-1
        \right)
        \alpha_m(\bfs)w_m^2 .
        \label{eq:delta_arrival_exact_bound}
    \end{align}
    Since $s_j=0$, the routing policy implies  $w_m < \tau_{mj}$.
    Also, since $\kappa_{m\ell}^*\in(0,1)$, $\alpha_m(\bfs) \leq \lambda_m$ by~\eqref{eq:def_test_function}.
    Consequently,
    \begin{align}
        0
        \leq
        \Delta_i^{\mathrm{arr}}(\bfx)
        &\leq
        \max_{j'\in\calJ}
        \Bigl(
            \sum_{m\in\calI:\,(mj')\in\calL_T}
            \bigl(
                \tfrac{1}{\kappa_{mj'}^*}-1
            \bigr)
            \lambda_m \tau_{mj'}^2
        \Bigr)
        <\infty .
        \label{eq:bnd_delta_arr}
    \end{align}
    Summing over all empty queues gives 
    \begin{align}
        \label{eq:bnd_change_arr}
        0
        \leq
        \sum_{i\in\calI:q_i=-}\lambda_i\Delta_i^{\mathrm{arr}}(\bfx)
        &\leq
        \max_{j'\in\calJ}
        \Bigl(
            \sum_{m\in\calI:\,(mj')\in\calL_T}
            \bigl(
                \tfrac{1}{\kappa_{mj'}^*}-1
            \bigr)
            \lambda_m \tau_{mj'}^2
        \Bigr)
        \Bigl(\sum_{i\in\calI} \lambda_i\Bigr)
        =:
        c_1
        <\infty .
    \end{align}

    We next analyze $\Delta_j^{\textrm{dep}}$. 
    Let $j\in\calJ$ with $s_j=1$.
    By~\eqref{eq:def_Q_int_departure},
    \begin{align}
        \label{eq:delta_dep_split}
        \Delta_j^{\textrm{dep}}
        &=
        \begin{cases}
            f(\bfq,\bfs^{j\mapsto 0},\bfw) - f(\bfq,\bfs,\bfw),
            & \calC_j(\bfx)=\emptyset 
            ,
            \\
            \int_E
            \bigl(f(\bfy)-f(\bfx)\bigr) c_{m_j(\bfx),j}(\bfx,\rmd \bfy),
            & \calC_j(\bfx)\neq\emptyset
            .
        \end{cases}
    \end{align}
    By~\eqref{eq:def_test_function} and~\eqref{eq:alpha_property},
    \begin{align}
        \label{eq:delta_dep_sum}
        f(\bfq,\bfs^{j\mapsto 0},\bfw) - f(\bfq,\bfs,\bfw)
        &=
        \sum_{i\in\calI}
        \bigl(
            \alpha_i(\bfs^{j\mapsto 0})
            -
            \alpha_i(\bfs)
        \bigr)
        w_i^2
        =
        \sum_{i\in\calI:\,(ij)\in\calL_T}
        \left(
            \kappa_{ij}^*-1
        \right)
        \alpha_i(\bfs)w_i^2
        \leq 0
        .
    \end{align}
    Here, the last inequality follows since $\kappa_{ij}^* \in (0,1)$ for all $(ij)\in\calL_T$, and since $\alpha_i(\bfs), w_i \geq 0$. 
    Hence, $\Delta_j^{\textrm{dep}} \leq 0$ if $\calC_j(\bfx) = \emptyset$. 
    Suppose that $\calC_j(\bfx)\neq \emptyset$ and write $m = m_j(\bfx)$, then by~\eqref{eq:def_c_line}, \eqref{eq:def_test_function}, and \eqref{eq:delta_dep_split}, using that the coefficient $\alpha$ only depend on $\bfs$, 
    \begin{align}
        \Delta_j^{\textrm{dep}}
        &= 
        \bigl(1-F_{A_m}(w_m)\bigr)
        \bigl(
            f\bigl(
                \bfq^{m\mapsto -},
                \bfs,
                \bfw^{m\mapsto 0}
            \bigr)
            -
            f(
                \bfq,
                \bfs,
                \bfw
            )
        \bigr) \\
        &\ 
        +
        \int_0^{w_m}
        \bigl(
            f\bigl(
                \bfq,
                \bfs,
                \bfw^{m\mapsto w_m-a}
            \bigr)
            -
            f(
                \bfq,
                \bfs,
                \bfw
            )
        \bigr)
        \rmd F_{A_m}(a) \\
        &= 
        \alpha_m(\bfs)
        \Bigl(
        \bigl(1-F_{A_m}(w_m)\bigr)
        \bigl(
            0 - w_m^2
        \bigr)
        +
        \int_0^{w_m}
        \bigl(
            (w_m-a)^2 - w_m^2
        \bigr)
        \rmd F_{A_m}(a) 
        \Bigr) \\
        &=
        \alpha_m(\bfs)
        \Bigl(
        - w_m^2
        +
        \int_0^{w_m}
        (w_m-a)^2
        \rmd F_{A_m}(a)
        \Bigr)
        .
    \end{align}
    Using $F_{A_m}(x) = 1 - \exp(-\lambda_m x)$ and integration by parts gives 
    \begin{align}
        \label{eq:bnd_change_dep}
        \Delta_j^{\textrm{dep}}
        &= 
        \alpha_m(\bfs)
        \Bigl(
            - \frac{2w_m}{\lambda_m}
            +
            \frac{2}{\lambda_m^2}
            \bigl(1-\exp(-\lambda_m w_m)\bigr)
        \Bigr)
        \leq 
        \alpha_m(\bfs)
        \Bigl(
            - \frac{2w_m}{\lambda_m}
            +
            \frac{2}{\lambda_m^2}
        \Bigr)
        ,
    \end{align}
    where the last inequality follows since $w_m, \lambda_m\geq 0$. 

    We obtain from~\eqref{eq:drift_decomposition_proof},~\eqref{eq:bnd_change_arr}, and~\eqref{eq:bnd_change_dep}, 
    \begin{align}
        \scrA_\bfX f(\bfx)
        &\leq
        \sum_{i\in\calI}  2 \alpha_i(\bfs)  w_i
        - 
        \sum_{j\in\calJ: s_j=1, \calC_j(\bfx)\neq\emptyset}
        \frac{2 \mu_j \alpha_{m_j(\bfx)}(\bfs)\, w_{m_j(\bfx)}}{\lambda_{m_j(\bfx)}}
        + \;
        c_2,
        \label{eq:assembled_bound_sq} 
    \end{align}
    where 
    $c_2 := c_1 +  \max_{m\in\calI} (2/ \lambda_{m}^2) \sum_{j\in\calJ}\mu_j$.

    Let $\tau_{\max}:=\max_{(ij)\in\calL}\tau_{ij}<\infty$. 
    We next show that
    \begin{align}
        \sum_{j\in\calJ: s_j=1,\ \calC_j(\bfx)\neq\emptyset}
        \frac{\alpha_{m_j(\bfx)}(\bfs)w_{m_j(\bfx)}}{\lambda_{m_j(\bfx)}}
        &\geq
        \sum_{j\in\calJ}
        \max_{i:(ij)\in\calL}(w_i-\tau_{\max})^+ .
        \label{eq:service_lower_sum}
    \end{align}
    Fix $j\in\calJ$ with $s_j=1$ and $\calC_j(\bfx)\neq\emptyset$, and write $m:=m_j(\bfx)$. 
    We argue that $\alpha_m(\bfs)w_m/\lambda_m\geq \max_{i:(ij)\in\calL}(w_i-\tau_{\max})^+$.
    If the maximum is zero, this is immediate. Otherwise, there exists $i$ with $(ij)\in\calL$ and $w_i>\tau_{\max}$. Since $\tau_{ij}\leq\tau_{\max}$, type~$i$ is eligible for server $j$, and hence $i\in\calC_j(\bfx)$. 
    Since $m=m_j(\bfx)$, we have $w_m=\max_{k\in\calC_j(\bfx)}w_k$, so every $i$ with $w_i > \tau_{\max}$ satisfies $w_i\leq w_m$. 
    Therefore, $\max_{i:(ij)\in\calL}(w_i-\tau_{\max})^+\leq w_m-\tau_{\max}$.
    Moreover, since $w_m>\tau_{\max}$, class $m$ has no idle compatible threshold server, so $\calJ_m^T(\bfs)=\emptyset$ and therefore $\alpha_m(\bfs)=\lambda_m$ by~\eqref{eq:def_test_function}. 
    Hence, $\alpha_m(\bfs)w_m/\lambda_m=w_m\geq w_m-\tau_{\max}\geq \max_{i:(ij)\in\calL}(w_i-\tau_{\max})^+$.

    It remains only to justify that the right side in \eqref{eq:service_lower_sum} may be summed over all $j\in\calJ$. 
    If $\calC_j(\bfx) = \emptyset$ for a server~$j$, then there are no eligible customers, hence $\max_{i:(ij)\in\calL}(w_i - \tau_{\max})^+ = 0$.
    Similarly,  $s_j = 0$ implies $\max_{i:(ij)\in\calL}(w_i - \tau_{\max})^+ = 0$, since server $j$ cannot be idle while an eligible compatible customer is waiting. 
    Hence, all servers with $\calC_j(\bfx) = \emptyset$ or $s_j = 0$ have zero contribution on the right side of \eqref{eq:service_lower_sum}.

    Since $\calI$ is finite and the inequalities in \eqref{eq:stab_cond_general} are strict, there exists $\eta>0$ such that for all $\calA \subseteq \calI$,
    \begin{align}
        \label{eq:hall_eta}
        (1+ \tfrac 12 \eta)
        \sum_{i\in \calA} \lambda_i
        \leq 
        \sum_{j\in\calS(\calA)} \mu_j
        .
    \end{align}
    For $t\geq 0$, define 
    $
        \calA_t
        :=
        \{i\in\calI:(w_i-\tau_{\max})^+>t\}
        ,
    $
    then $\max_{i:(ij)\in\calL} (w_i-\tau_{\max})^+ > t$ if and only if $j\in\calS(\calA_t)$.
    Using this together with $a = \int_0^\infty \ind{a>t} \, \rmd t$ for any $a\geq 0$, and finiteness of the summations gives
    \begin{align}
        \sum_{j\in\calJ}
        \mu_j
        \max_{i:(ij)\in\calL}
        (w_i-\tau_{\max})^+ 
        &=
        \sum_{j\in\calJ}
        \mu_j
        \int_0^\infty
        \ind{
            \max_{i:(ij)\in\calL}
            (w_i-\tau_{\max})^+ > t
        }
        \,\rmd t
        \notag \\
        &=
        \int_0^\infty
        \sum_{j\in\calJ}
        \mu_j
        \ind{
            \max_{i:(ij)\in\calL}
            (w_i-\tau_{\max})^+ > t
        }
        \,\rmd t
        \notag \\
        &=
        \int_0^\infty
        \sum_{j\in\calS(\calA_t)}\mu_j
        \,\rmd t \\
        &\stackrel{\eqref{eq:hall_eta}}{\geq}
        (1+ \tfrac 12 \eta)
        \int_0^\infty
        \sum_{i\in \calA_t}\lambda_i
        \,\rmd t \\
        &=
        (1+ \tfrac 12 \eta)
        \sum_{i\in\calI}
        \lambda_i (w_i-\tau_{\max})^+ 
        .
        \label{eq:bnd_level_set_mid}
    \end{align}
    Since
    $
        (w_i-\tau_{\max})^+
        \geq
        w_i-\tau_{\max}
        ,
    $
    we obtain from~\eqref{eq:bnd_level_set_mid}
    \begin{align}
        \label{eq:hall_level_bound}
        \sum_{j\in\calJ}
        \mu_j
        \max_{i:(ij)\in\calL}
        (w_i-\tau_{\max})^+
        \geq
        (1+ \tfrac 12 \eta)
        \sum_{i\in\calI}\lambda_i w_i
        -
        (1+\tfrac 12 \eta)\tau_{\max}\sum_{i\in\calI}\lambda_i .
    \end{align}
    Lastly, since $\kappa_{ij}^* \in (0,1)$ we have $\lambda_i \geq \alpha_i(\bfs) $ for all $i\in\calI$ by~\eqref{eq:def_test_function}.
    We thus find by \eqref{eq:service_lower_sum}  and~\eqref{eq:hall_level_bound}
    \begin{align}
        \sum_{j\in\calJ: s_j=1, \calC_j(\bfx)\neq\emptyset}
        \frac{\mu_j\alpha_{m_j(\bfx)}(\bfs)w_{m_j(\bfx)}}
            {\lambda_{m_j(\bfx)}}
        &\geq
        (1+ \tfrac 12 \eta)
        \sum_{i\in\calI}\lambda_i w_i
        -
        c_3
        \geq
        (1+ \tfrac 12 \eta)
        \sum_{i\in\calI} \alpha_i(\bfs)w_i
        -
        c_3,
        \label{eq:bnd_dep_eta}
    \end{align}
    with
    $
        c_3
        :=
        (1+ \tfrac 12 \eta)\tau_{\max}\sum_{i\in\calI}\lambda_i .
    $
    The proof is concluded from~\eqref{eq:assembled_bound_sq} and~\eqref{eq:bnd_dep_eta}, with $c_0 = c_2+c_3$. 

\end{proof}

\subsection{Proof of \Cref{lem:petite}}
\label{app:proof_lem_petite}

\begin{proof}
    Let $a\geq 0$ and let $\bfx_0$ denote the all-empty state, i.e., 
    $
        \bfx_0:=(\bfq_0,\bfzero,\bfzero),
    $
    with $\bfq_{0,i}=-$ for all $i\in\calI$. 
    Let
    $
        P^t(\bfx,A):=\P_\bfx(\bfX(t)\in A)
    $
    denote the transition kernel of $\bfX(t)$.
    To show that $B_a$ is small, it suffices to show that there exist $t>0$ and $\eps_a>0$ such that
    $
        P^t(\bfx,\cdot)
        \geq
        \eps_a \delta_{\bfx_0}(\cdot)
    $
    for all $\bfx\in B_a$, where $\delta_{\bfx_0}(\cdot)$ is the Dirac measure at $\bfx_0$. 

    Let
    $
        \tau_{\max}:=\max_{(ij)\in\calL}\tau_{ij}<\infty,
    $
    $
        \lambda_{\Sigma}:=\sum_{i\in\calI}\lambda_i,
    $
    $
        \mu_{\min}:=\min_{j\in\calJ}\mu_j>0,
    $
    and
    $ 
        t_0:=(|\calI| + |\calJ|)(\tau_{\max}+2)+1 
        .
    $

    We construct an event on which the system is empty at time $t_0$.
    First require that there are no arrivals during $[0,t_0]$. 
    This occurs with probability
    $
        \exp(-\lambda_{\Sigma}t_0)>0.
    $
    Conditional on the event that there are no arrivals in $[0,t_0]$, the only customers that can be present in the system at time $t_0$ are those already in service  and those waiting in the non-empty queues at time $0$.

    By construction, $t_0$ consist of $|\calI|+|\calJ|$ blocks of length $\tau_{\max}+2$.
    We call a block successful if either a busy server becomes idle, or a non-empty queue becomes empty during that block.

    If at least one server is busy at the start of a block, then a departure occurs within the next time unit with probability at least
    $
        1-\exp(-\mu_{\min})>0.
    $
    After such a departure, either the server becomes idle, or, if an eligible customer is present, starts a new service.
    Since all thresholds are bounded by $\tau_{\max}$, within time at most $\tau_{\max}$ some non-empty queue becomes eligible for a compatible idle server. 
    The probability that its service time is at most one time unit is at least
    $
        1-\exp(-\mu_{\min}) 
        .
    $
 
    Suppose a type-$i$ \gls{FIL} customer with waiting time $w_i$ is taken into service. 
    Then queue~$i$ becomes empty if  the interarrival time $A_i \sim \expdist{\lambda_i}$ exceeds $w_i$. 
    During the interval $[0,t_0]$, all waiting times are bounded by $a+t_0$ for $\bfx\in B_a$. 
    Hence, uniformly over all states reached before time $t_0$,
    \begin{align}
        \P(A_i>w_i)
        \geq
        \P(A_i>a+t_0)
        =
        \exp(-\lambda_i(a+t_0))
        \geq
        \exp(-\lambda_{\Sigma}(a+t_0))
        > 0
        .
    \end{align}
    Thus there is a uniformly positive probability that the queue becomes empty.

    Each successful block reduces the number of initially present customers that still have to be completed or started into service. 
    Initially, there are at most $|\calJ|$ customers in service and at most $|\calI|$ non-empty queues.
    Therefore, after at most
    $
        |\calI|+|\calJ|
    $
    successful blocks, all initially present work has been removed and all queues
    are empty. 
    Since no arrivals occur on $[0,t_0]$, the process then stays in $\bfx_0$ until time $t_0$.

    Combining the uniform lower bounds above, we have
    $
        \P_{\bfx}(\bfX(t_0)=\bfx_0)
        \geq
        \eps_a
        ,
    $
    for all $\bfx\in B_a$,
    where
    \begin{align}
        \eps_a
        :=
        \exp(-\lambda_{\Sigma}t_0)
        (1-\exp(-\mu_{\min}))^{|\calJ|}
        \exp(-|\calI|\lambda_{\Sigma}(a+t_0))
        > 0
        .
    \end{align}
    Hence, for every $A\in\scrB(E)$,
    \begin{align}
        P^{t_0}(\bfx,A)
        \geq
        \P_{\bfx}(\bfX(t_0)=\bfx_0)\delta_{\bfx_0}(A)
        \geq
        \eps_a\delta_{\bfx_0}(A),
        \qquad 
        \bfx\in B_a .
    \end{align}
    This concludes the proof that $B_a$ is small.
\end{proof}

\end{document}